\documentclass[10pt]{ifacconf}

\usepackage{natbib}            
\usepackage{graphicx}          
\usepackage{graphics} 
\usepackage{amsmath} 
\usepackage{amssymb}  
\usepackage{eurosym}

\newtheorem{remark}{Remark}
\newtheorem{problem}{Problem}

\newtheorem{proof}{Proof}

\begin{document}

\begin{frontmatter}

\title{Some remarks on \\ wheeled autonomous vehicles \\ and the evolution of their control design\thanksref{footnoteinfo}} 

\thanks[footnoteinfo]{Work partially supported by the French national project INOVE/ANR 2010 BLANC 308.}

\author[Premier]{Brigitte d'ANDR\'EA-NOVEL}
\author[Deuxieme]{Lghani MENHOUR}
\author[Troisieme,Cinquieme]{Michel FLIESS}
\author[Quatrieme]{Hugues MOUNIER}

\address[Premier]{Centre de Robotique, MINES ParisTech \& PSL Research University, 60 boulevard Saint-Michel, 75272 Paris cedex 06, France. \\ (e-mail: brigitte.dandrea-novel@mines-paristech.fr)}
\address[Deuxieme]{Universit\'e de Reims Champagne-Ardenne, IUT de Troyes, 9 rue du Qu\'ebec, 10026 Troyes, France. (e-mail: lghani.menhour@univ-reims.fr)} 
\address[Troisieme]{LIX (CNRS, UMR 7161), \'Ecole polytechnique, 91128 Palaiseau, France. (e-mail: Michel.Fliess@polytechnique.edu)}
\address[Quatrieme]{L2S (UMR 8506), CNRS -- Sup\'elec -- Universit\'{e} Paris-Sud, 3 rue Joliot-Curie, 91192 Gif-sur-Yvette, France. \\ (e-mail: hugues.mounier@lss.supelec.fr)}
\address[Cinquieme]{AL.I.E.N. (ALg\`{e}bre pour Identification et Estimation Num\'{e}riques), 24-30 rue Lionnois, BP 60120, 54003 Nancy, France.\\ (e-mail: michel.fliess@alien-sas.com)}


\begin{abstract}                          
Recent investigations on the longitudinal and lateral control of wheeled autonomous vehicles are reported. Flatness-based techniques are first introduced via a simplified model. It depends on some physical parameters, like cornering stiffness coefficients of the tires, friction coefficient of the road, \dots, which are notoriously difficult to identify. Then a model-free control strategy, which exploits the flat outputs, is proposed. Those outputs also depend on physical parameters which are poorly known, \textit{i.e.}, the vehicle mass and inertia and the position of the center of gravity. A totally model-free control law is therefore adopted. It employs natural output variables, namely the longitudinal velocity and the lateral deviation of the vehicle. This last method, which is easily understandable and implementable, ensures a robust trajectory tracking problem in both longitudinal and lateral dynamics. Several convincing computer simulations are displayed.
\end{abstract}

\end{frontmatter}

\section{Introduction}
 The lateral and longitudinal control of wheeled autonomous vehicles is an important topic which has already attracted many promising studies (see, \textit{e.g.}, \cite{Ackermann95, Novel01,anton,att,Choi09, Cerone09, Chou05, Fuchsumer2005, Marino09, Martinez07, cep14, Nouveliere02, Poussot11a, Rajamani00, Villagra09, Villagra11, vil2, Zheng06}, \dots), where various advanced theoretical tools are utilized. This short communication does not permit unfortunately to summarize them. Let us nevertheless notice that most of them are model-based. The aim of this presentation is to explain and justify the evolution of our viewpoint which started with a flatness-based setting (\cite{cep14}), \textit{i.e.}, a model-based approach. It is now adopting a fully model-free standpoint (see \cite{Menhour13b} and \cite{ECC2015}). As a matter of fact severe difficulties are encountered to
\begin{itemize}
\item write a mathematical model which takes into account all the numerous complex phenomena,
\item calibrate the existing model in various changing situations like cornering stiffness coefficients of the tires, and friction coefficient of the road. 
\end{itemize}
Our paper is organized as follows. Section \ref{Section_2} presents a nonlinear longitudinal and lateral flatness-based control. In Section \ref{Section_3} the two different model-free control strategies are developed.
Simulation results with noisy data and suitable reference trajectories acquired on a track race, are displayed in Section \ref{Section_4}. Let us emphasize that the second model-free control is quite robust.  Concluding remarks may be found in Section \ref{Section_6}.
\section{Longitudinal and lateral flatness-based control}
\label{Section_2}
\subsection{3DoF NonLinear Two Wheels Vehicle Model}
\label{Section_2_1}

\begin{figure}[!ht]
\centering
\rotatebox{0}{
\includegraphics[scale= 0.43]{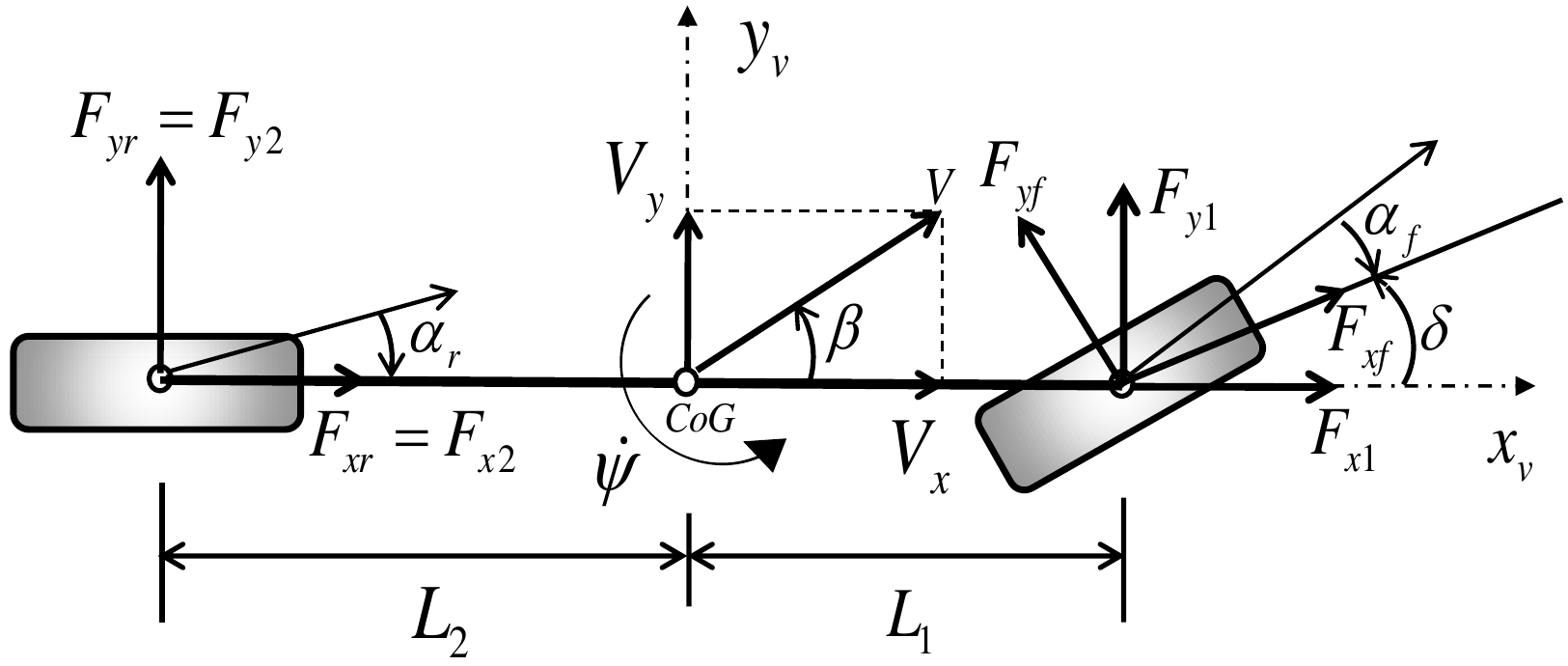}}
\caption{Nonlinear two wheels vehicle control model} \label{NLBM1}
\end{figure}

The 3DoF-NLTWVM in Figure \ref{NLBM1}, which is used to design the combined control law, provides an interesting approximation of the longitudinal and lateral dynamics of the vehicle in normal driving situations. See Table \ref{notations_vehicle}
for the notations.
{\scriptsize
\begin{table}[!ht]
\caption{}
\label{notations_vehicle}
\centering
\begin{tabular}{|l|l|}%
\hline
   Symbol & Variable name \\
   \hline
   $V_x$ & longitudinal speed [$km.h$]\\
   $V_y$ & lateral speed [$km.h$]\\
   $a_x$ & longitudinal acceleration [$m/s^2$]\\
   $a_y$ & lateral acceleration [$m/s^2$]\\
   $\dot{\psi}$& yaw rate $[rad/s]$\\
   $\psi$ & yaw angle $[rad]$\\
   $\beta$ & sideslip angle $[rad]$\\
   $\alpha_{f,\,r}$ & front and rear tire slip angles $[rad]$\\
   $\omega_i$ & wheel angular speed of the wheel $i$ $[rad/s]$\\
   $T_{\omega}$ & wheel torque $[Nm]$\\
   $\delta $ & wheel steer angle $[deg]$\\
   $C_f$, $C_r$ & front and rear cornering stiffnesses $[N.rad^{-1}]$\\
   $F_{xi}$ & longitudinal forces in vehicle coordinate $[N]$\\
   $F_{yi}$ & lateral forces in vehicle coordinate $[N]$\\
   $F_{xf}$ & front longitudinal forces in tire coordinate $[N]$\\
   $F_{yf}$ & front lateral forces in tire coordinate $[N]$\\
   $R$ & tire radius $[m]$ \\
   $g$ & acceleration due to gravity $[m/s^2]$\\
   $L_f$ & distances from the CoG to the front axles $[m]$\\
   $L_r$ & distances from the CoG to the rear axles $[m]$\\
   $I_z$ & yaw moment of inertia $[kg m^2]$\\
   $I_r$ & wheel moment of inertia $[kg m^2]$\\
   $m$ & vehicle mass $[kg m^2]$\\
   $M_z$ & yaw moment $[Nm]$\\
   \hline
\end{tabular}
\end{table}
}

\noindent The corresponding dynamical equations read:
{\small
\begin{equation*}
\label{Non_linear_bicycle_model}
\left \{ \begin{array}{l}
m a_x=m(\dot{V}_x-\dot{\psi}V_y)= (F_{x1}+F_{x2}) \\
m a_y=m(\dot{V}_y+\dot{\psi}V_x)=(F_{y1}+F_{y2}) \\
I_z\ddot{\psi}= {M_{z1}+M_{z2}}
\end{array}
\right.
\end{equation*}
}
\noindent If a linear tire model with small slip angles is assumed (see \cite{cep14} for more details), $C_r$ (resp. $C_f$) denoting the cornering stiffness coefficient of the rear (resp. front) wheel, and
$\omega _r$ (resp. $\omega_f$ ) being the angular velocity of the rear (resp. front) wheel, the previous system can be rewritten:
\begin{equation}
\label{affine_NL_modele}
\dot{x}=f(x,t)+g(x,t)u
\end{equation}
with
$$
\scriptsize{\label{f_x_t}
f(x,t)=\left[
\begin{array}{c}
\dot{\psi}V_y  - \frac{I_r}{mR}(\dot{\omega}_r + \dot{\omega}_f)  \\[2mm]
-\dot{\psi}V_x +\frac{1}{m}\left( -C_f\left(\frac{V_y+L_f\dot{\psi}}{V_x}\right) - C_r\left(\frac{V_y-L_r\dot{\psi}}{V_x}\right)\right) \\[2mm]
\frac{1}{I_z}\left(-L_fC_f\left(\frac{V_y+L_f\dot{\psi}}{V_x}\right)+L_r C_r\left(\frac{V_y-L_r\dot{\psi}}{V_x}\right)\right)
\end{array}
\right]}
$$
$$
g(x,t) = \left[
\begin{array}{cc}
\frac{1}{mR} & \frac{C_f}{m}\left(\frac{V_y+L_f\dot{\psi}}{V_x}\right)\\ [1mm]
0 & ({C_fR-I_r\dot{\omega}_f})/{mR}\ \\ [1mm]
0 & ({L_fC_f R-L_fI_r\dot{\omega}_f})/{I_zR}
\end{array}
\right], \hspace{0.1cm}
x=
\left[
\begin{array}{c}
V_x\\ V_y \\ \dot{\psi}
\end{array}
\right]
$$ 
$u= \left[ \begin{array}{cc} u_1=T_{\omega},& u_2=\delta \end{array} \right]^T$. 
\subsection{Flatness property}
\label{Flatness_property_definition}
The system $\dot{x}=f(x,u)$, where $x=(x,\cdots,x_n)\in \mathbb{R}^{n}$ and $u=(u,\cdots,u_m)\in
\mathbb{R}^{m}$, is said to be \emph{differentially flat} (see \cite{Fliess95,Fliess99}, and \cite{murray,levine,hsr}) if, and only if,
\begin{itemize}
\item there exists a vector-valued function $h$ such that
\begin{equation}
\label{flat_outputs_general}
y=h(x,u,\dot{u},\cdots,u^{(r)})
\end{equation}
where $y=(y,\cdots,y_m)\in \mathbb{R}^{m}$, $r\in \mathbb{N}$;
\item the components of $x=(x,\cdots,x_n)$ and
$u=(u,\cdots,u_m)$ may be expressed as
\begin{equation}
\label{x_y_flat_sys}
x=A(y,\dot{y},\cdots,y^{(r_x)}),\>\>\> r_x\in \mathbb{N}
\end{equation}
\begin{equation}
\label{u_y_flat_sys}
u=B(y,\dot{y},\cdots,y^{(r_u)}),\>\>\> r_u\in \mathbb{N}
\end{equation}
\end{itemize}
Remember that $y$ in Equation \eqref{flat_outputs_general} is called a \emph{flat output}.
\subsection{Flatness-based longitudinal and lateral control}
\label{Section_2_2}
\begin{problem}
Introduce the outputs: 
\begin{equation}
\label{flatness_outputs}
\left \{\begin{array}{l}
y_1=V_x \\ [1mm]
y_2= L_f m V_y -I_z \dot{\psi}
\end{array}
\right.
\end{equation}
We want to show that the longitudinal speed $y_1$ and  the angular momentum $y_2$ of a point on the axis between the centers of the front and rear axles are flat outputs. \end{problem}
\begin{proof}
Some algebraic manipulations (see \cite{cep14} for more details) yield:
\begin{equation}
\label{x_A_y1_y2_y2p}
\begin{array}{c}
x=
\left[
\begin{array}{ccc}
V_x & V_y & \dot{\psi}
\end{array}
\right]^T
\begin{array}{c}
=A(y_1,y_2,\dot{y}_2)=
\end{array}\\
\small{
\left[
\begin{array}{c}
y_1\\ [2mm]
\frac{y_2}{L_f m} -  \frac{I_z }{L_f m} \left(\frac{L_f m y_1\dot{y}_2 + C_r(L_f+L_r)y_2 }{C_r(L_f+L_r)(I_z -L_rL_fm)+ (L_f my_1)^2 }\right)
\\ [2mm]
-\left(\frac{L_f m y_1\dot{y}_2 + C_r(L_f+L_r)y_2 }{C_r(L_f+L_r)(I_z -L_rL_fm)+ (L_f my_1)^2 }\right)
\end{array}
\right]}
\end{array}
\end{equation}
and
\begin{equation}
\label{U_B_y1_y2_y2p}
\begin{array}{c}
\left[
\begin{array}{c}
\dot{y}_1\\
\ddot{y}_2
\end{array}
\right]
= \Delta(y_1,y_2,\dot{y}_2) \left(
\begin{array}{c}
u_1\\
u_2
\end{array}
\right)+\Phi(y_1,y_2,\dot{y}_2) 
\end{array}
\end{equation}
The flatness property holds if the matrix $\Delta(y_1,y_2,\dot{y}_2)$ is invertible:
\begin{equation}
\label{delta_matrice}
\scriptsize{
\begin{array}{c}
\det (\Delta (y_1,y_2,\dot{y}_2) ) = \Delta_{11} \Delta_{22}- \Delta_{21} \Delta_{12} = \\[2mm]
\frac{\left(I_\omega \dot{\omega}_f - C_fR\right)\left(L_f^2y_1^2m^2 - Cr(L_f+L_r)L_rL_fm + C_rI_zL\right)}{I_zR^2 y_1 m^2}\neq 0
\end{array}}
\end{equation}
This determinant, which only depends on the longitudinal speed
$y_1=V_x$, is indeed nonzero:
\begin{itemize}
\item Acceleration of the wheel rotation is less than ${RC_f}/{I_\omega} \approx 10^4$. Thus
\begin{equation}
\label{cond_1}
I_\omega \dot{\omega}_f - C_f R\neq 0
\end{equation}
\item From $I_z > L_fm$, we deduce that
\begin{equation}
\label{cond_2}
C_r(L_r+L_f)(I_z - L_fm) + L_f^2m^2y_1^2 \neq 0 
\end{equation}
\end{itemize}
Thus
\begin{equation}
\label{u_B_y1_y1p_y2_y2p_y2pp}
\begin{array}{c}
u=\left[
\begin{array}{c}
T_\omega\\
\delta
\end{array}
\right]=
\begin{array}{c}
B(y_1,\dot{y}_1,y_2,\dot{y}_2, \ddot{y}_2) =
\end{array}\\[2mm]
\Delta^{-1}(y_1,y_2,\dot{y}_2)
\left(
\left[
\begin{array}{c}
\dot{y}_1\\
\ddot{y}_2
\end{array}
\right]
- \Phi(y_1,y_2,\dot{y}_2)
\right)
\end{array}
\end{equation}
with $r_x=1$ and $r_u=2$.
\end{proof}
\subsubsection{A tracking feedback control}
In order to track the desired output trajectories $y_1^{ref}$ and $y_2^{ref}$, set
    \begin{equation}
\label{lin_contro}
  \left [
  \begin{array}{c}
    \dot{y}_1 \\
    \ddot{y}_2
  \end{array} \right]=
  \left [
  \begin{array}{c}
  \dot{y}_1^{ref}+K_1^1e_{y_1}+K_1^2\int{e_{y_1} dt}\\[2mm]
\ddot{y}_2^{ref}+K_2^1\dot{e}_{y_2}+K_2^2e_{y_2}+K_2^3\int{e_{y_2} dt}
  \end{array} \right]
    \end{equation}
where, $e_{y_1}=y_1^{ref}-y_1=V_x^{ref}-V_x$ and
$e_{y_2}=y_2^{ref}-y_2$. Choosing the gains  $K_1^1$,
$K_1^2$, $K_2^1$, $K_2^2$ and $K_2^3$ is straightforward.
\subsection{Algebraic nonlinear estimation}
The control law contains derivatives of measured signals, which are of course noisy. This estimation is performed 
using the recent advances in \cite{Fliess08, Mboup09,sira-est}. It yields the following formulae (\cite{gretsi}), and, therefore, simple linear digital filters:
\begin{itemize}
    \item Denoising:
\begin{equation}
\label{filtering}
\hat{y}(t)=\frac{2!}{T^2} \int^{t}_{t-T}(3(t-\tau)-T)y(\tau)d\tau
\end{equation}
\item Numerical differentiation of a noisy signal:
\begin{equation}
\label{Derivative}
\hat{\dot{y}}(t)=-\frac{3!}{T^3} \int^{t}_{t-T}(2T(t-\tau)-T)y(\tau)d\tau
\end{equation}
\end{itemize}
The sliding time window $\left[t-T,\>t \right]$ may be quite short.
\section{Model-Free design for vehicle control}
\label{Section_3}
\subsection{Model-free design: motivation}
\begin{figure}[!ht]
\centering
\includegraphics[scale=0.5]{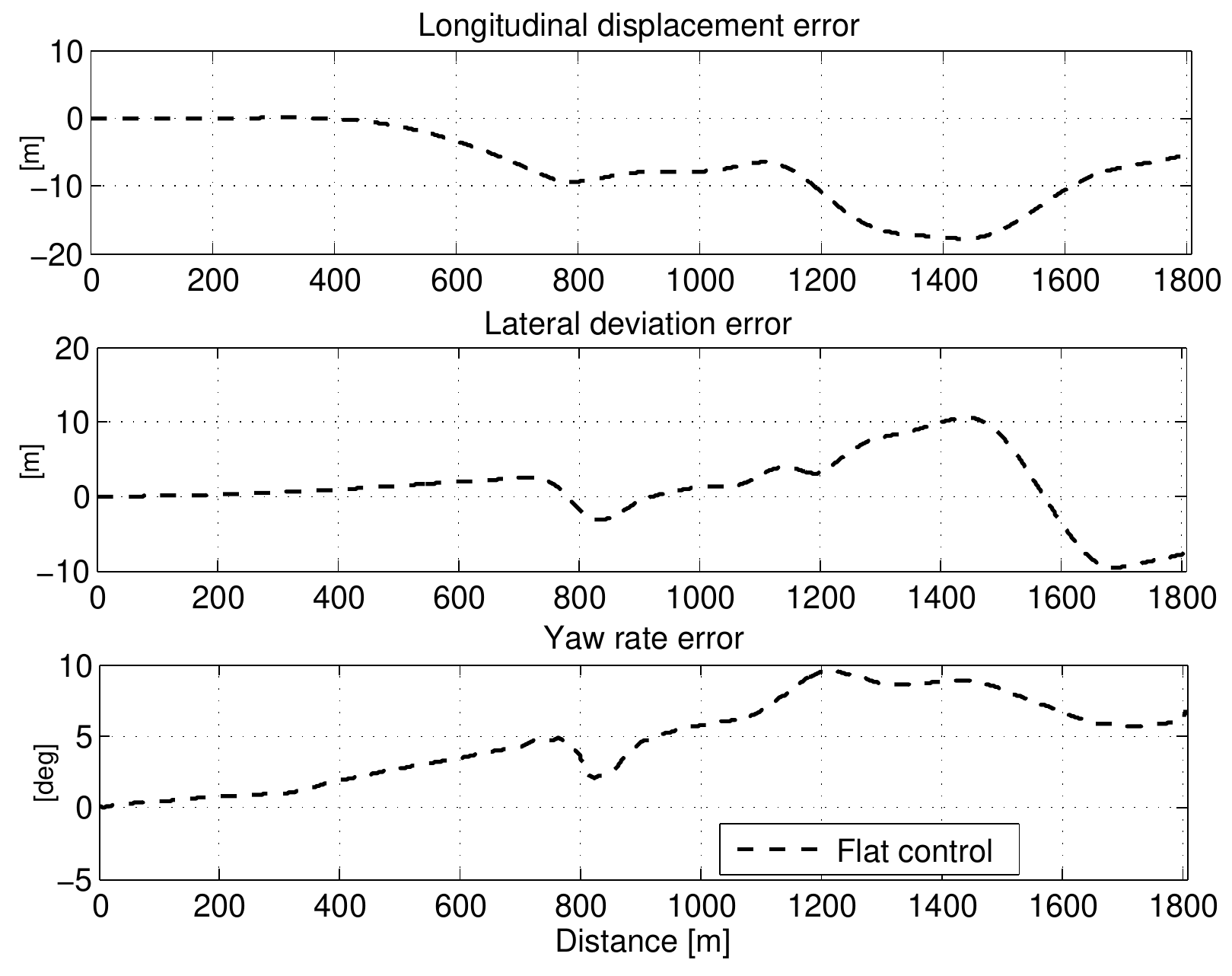}
\vspace{-0.4cm}
\caption{Simulation test with $0.3C_f$ and $0.3Cr$: trajectory tracking with nonlinear flatness-based vehicle control}
\label{Test_flat_control_70_Cf_Cr_V1}
\end{figure}

\begin{figure}[!ht]
\centering
\includegraphics[scale= 0.5]{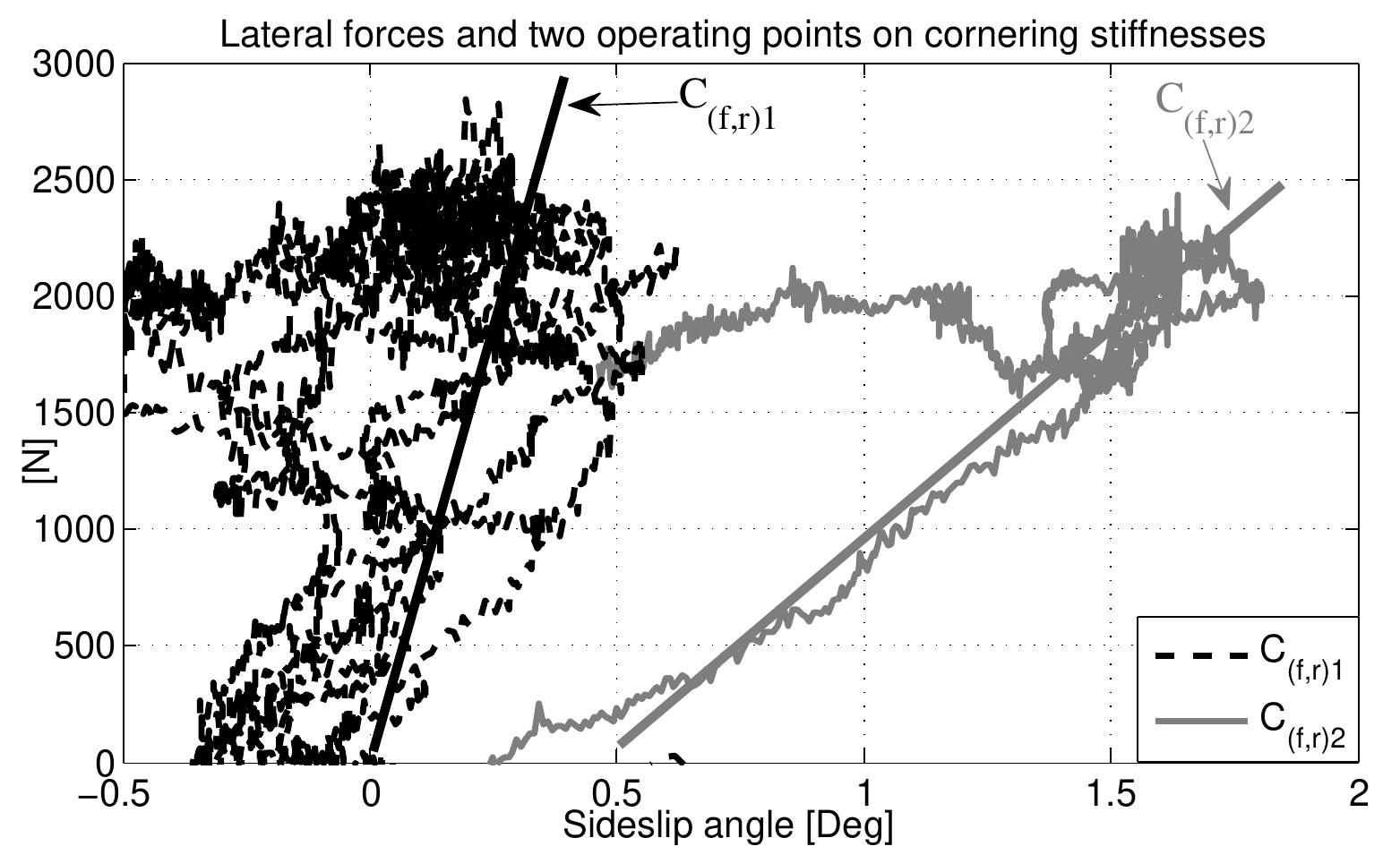}
\vspace{-0.4cm}
\caption{Experimental lateral force characteristic and variations of cornering stiffness coefficients ($C_f$ and $C_r$)}
\label{switching_strategy_on_Cf_Cr}
\end{figure}

Realistic vehicle models are usually unsuitable for an implementable control design. The above nonlinear flat controller provides good results with known nominal values of the cornering stiffnesses $C_f$ and $C_r$. Uncertainties on $C_f$ and $C_r$ yield unsatisfactory closed-loop performances: see Figure \ref{Test_flat_control_70_Cf_Cr_V1}.
It is well known that $C_f$ and $C_r$ are very sensitive to high dynamic loads. Figure \ref{switching_strategy_on_Cf_Cr} displays different dynamics of the front and rear cornering stiffness coefficients $C_f$ and $C_r$ during an experimental braking maneuver.
\subsection{A short review on model-free control\protect\footnote{See \cite{ijc13} for more details and further explanations. See \cite{nice} for a cheap hardware implementation. Both references show the online 
character of this setting, This fact has been confirmed by many concrete case-studies.}}
\label{Section_3_1}
Model-free control (\cite{ijc13}) was already successfully used in many concrete case-studies (see, \textit{e.g.}, \cite{ijc13}, and \cite{dk,toulon}; the references therein provide numerous other examples). Let us insist here on the already existing applications to vehicles: \cite{Choi09,formentin,ECC2015,Villagra09,vil2}. 
\subsubsection{The ultra-local model}\label{A}
Replace the unknown SISO system by the \emph{ultra-local model}:
\begin{equation}
y^{(\nu)} = F + \alpha u
\label{ultralocal}
\end{equation}
where 
\begin{itemize}
\item $\nu \geq 1$ is the derivation order,
\item $\alpha \in \mathbb{R}$ is chosen such that $\alpha u$ and $y^{(\nu)}$ are of the same order of magnitude, 
\item $\nu$ and $\alpha$ are chosen by the practitioner. 
\end{itemize}
\begin{remark}
In all the existing concrete examples $\nu = \ 1 \ \text{or} \ 2$. 
Magnetic bearings (\cite{Miras13}) with their low friction provide the only instance where $\nu = 2$ (see \cite{ijc13} for an explanation).
\end{remark}
Some comments on $F$ are in order:
\begin{itemize}
\item $F$ is estimated via the measurements of the input $u$ and the output $y$,
\item $F$ does not distinguish between model mismatches and perturbations.
\end{itemize}
\subsubsection{Intelligent controllers}
Set $\nu = 2$ in Equation \eqref{ultralocal}:
\begin{equation}
\label{MFC_4_n_2}
\ddot{y} = F+ \alpha u 
\end{equation}
The corresponding \emph{iPID}, \textit{i.e}, \emph{intelligent Proportional-Integral}   \emph{-Derivative controller}, reads
\begin{equation}
\label{iPID_c}
u = - \frac{\left( F - \ddot{y}^d + K_P e + K_I\int e dt + K_D \dot{e}\right)}{\alpha}  
\end{equation}
where
\begin{itemize}	
\item $y^d$ is the reference trajectory,
\item $e = y-y^d$ is the tracking error and $y^d$ is a desired signal, 
\item $K_P$, $K_I$, $K_D \in \mathbb{R}$ are the usual gains.
\end{itemize}
Combining Equations \eqref{MFC_4_n_2} and \eqref{iPID_c} yields
$$
\ddot{e} + K_P e + K_I \int e dt + K_D \dot{e}= 0
$$
where $F$ does not appear anymore. Gain tuning becomes therefore quite straightforward. This is a major benefit when compared to ``classic'' PIDs (see, \textit{e.g.}, \cite{astrom,murray}). 
If $K_I = 0$ we obtain the \emph{intelligent Proportional-Derivative controller}, or \emph{iPD},
\begin{equation}
\label{iPD_c}
u = - \frac{\left( F - \ddot{y}^d + K_P e dt + K_D \dot{e}\right)}{\alpha}  
\end{equation}
Set $\nu = 1$ in Equation \eqref{ultralocal}:
\begin{equation}
\label{MFC_4_n_1}
\dot{y} = F+ \alpha u 
\end{equation}
The corresponding \emph{intelligent Proportional-Integral controller}, or \emph{iPI}, reads: 
\begin{equation}
\label{iPI_c}
u = - \frac{\left(F - \dot{y}^d + K_P e + K_I\int e dt \right)}{\alpha}
\end{equation}	
If $K_I = 0$ in Equation \eqref{iPI_c}, we obtain the \emph{intelligent proportional controller}, or \emph{iP}, which, until now, turns out to be the most useful intelligent controller:
\begin{equation}
\label{ip}
{u = - \frac{F - \dot{y}^\ast + K_P e}{\alpha}}
\end{equation} 
\subsubsection{Algebraic estimation of $F$}
\label{F}
$F$ in Equation \eqref{ultralocal} is assumed to be ``well'' approximated by a piecewise constant function $F_{\text{est}} $. According to the algebraic parameter identification due to \cite{sira1,sira2} (see also \cite{sira-est}), rewrite, if $\nu = 1$, Equation \eqref{MFC_4_n_1} in the operational domain (see, \textit{e.g.}, \cite{Yosida84}) 
$$
s Y = \frac{\Phi}{s}+\alpha U +y(0)
$$
where $\Phi$ is a constant such that $\frac{\Phi}{s}$ is the operational transform of $F$ which is supposed to be constant on the sliding window. We get rid of the initial condition $y(0)$ by multiplying both sides on the left by $\frac{d}{ds}$:
$$
Y + s\frac{dY}{ds}=-\frac{\Phi}{s^2}+\alpha \frac{dU}{ds}
$$
Noise attenuation is achieved by multiplying both sides on the left by $s^{-2}$. It yields in the time domain the real time estimation
\begin{equation}
\label{integral}
 F_{\text{est}}(t)  =-\frac{6}{\tau^3}\int_{t-\tau}^t \left\lbrack (\tau -2\sigma)y(\sigma)+\alpha\sigma(\tau -\sigma)u(\sigma) \right \rbrack d\sigma 
\end{equation}
where $\tau > 0$ may be quite small. This integral may of course be replaced in practice by a classic digital filter. The extension to $\nu = 2$ is given by the following estimator:
\begin{equation}
\label{integral_2}
F_{\text{est}}(t) = 
  \left [
  \begin{array}{c}
-\frac{60}{\tau^5}\int_{t-\tau}^t(\tau^2 +6\sigma^2 - 6\tau \sigma)y_2(\sigma)d\sigma \\[2mm]
-\frac{30\alpha}{\tau^5} \int_{t-\tau}^t(\tau - \sigma)^2\sigma^2 u_2(\sigma)d\sigma
  \end{array} \right]
\end{equation}
\subsection{Model-free vehicle control and flat outputs}\label{mfcflat}
The only parameters supposed to be known are those used to compute $y_2$ in Equation \eqref{flatness_outputs}: the vehicle mass, the yaw moment of inertia and the position of the center of gravity. The control variables are: 
\begin{itemize}
\item the braking and traction wheel torques to control longitudinal motion,
\item the steering angle to control lateral and yaw motions. 
\end{itemize} 
The control variables $u_1=T_\omega$ and $u_2=\delta$ are respectively the braking and traction wheel torques, and the steering angle. From \eqref{u_B_y1_y1p_y2_y2p_y2pp}, the following two ultra-local models of longitudinal and lateral motions are derived: 
\begin{eqnarray*}
\label{MFC_LLVC2_1}
\text{longitudinal ultra-local model:} & \dot{y}_1 = F_1+ \alpha_1 u_1\\
\label{MFC_LLVC2_2}
\text{lateral ultra-local model:} & \dot{y}_2 = F_2+ \alpha_2 u_2 
\end{eqnarray*}
$u_1$ and $u_2$ represent respectively the wheel torque control of longitudinal motion $y_1$ and the steering angle control of lateral and yaw motions $y_2$. The tracking is achieved by two decoupled multivariable iPIs. Set, according to Equation \eqref{iPI_c}, for the
longitudinal iPI controller $u_1 = \frac{1}{\alpha_1} \left( F_1 + {\dot{y}}_1^{d} - K_P^{y_1} e_{y_1} - K_I^{y_1} \int e_{y_1} dt \right)$, and for the lateral iPI controller
$u_2 = \frac{1}{\alpha_2}\left( F_2 + {\dot{y}}_2^{d} - K_P^{y_2}e_{y_2} - K_I^{y_2}\int e_{y_2} dt \right)$,
where $e_{y_1}=y_1^{d}-y_1=V_x^{d}-V_x$ and $e_{y_2}=y_2^{d}-y_2$. The choice of the parameters $\alpha_1$, $\alpha_2$,  $K_P^{y_1}$, $K_I^{y_1}$, $K_P^{y_2}$ and $K_I^{y_2}$ is straightforward. Equations \eqref{integral} and \eqref{integral_2} yield respectively the estimation of $F_1$ and $F_2$.
\subsection{Model-free vehicle control and natural outputs}
In order to be totally parameter-independent, consider the following natural outputs, which can be obtained through direct measurements: 
 $$
    \left \{\begin{array}{l}
y_1= \text{longitudinal speed}\\ 
y_2= \text{lateral deviation}
\end{array}
\right.
  $$
The first output remains the measured (or well estimated) longitudinal velocity, and the second one is now given by the lateral deviation which is 
also supposed to be accessible from on line measurements. The two control variables are the same as in Section \ref{mfcflat}.
Newton's second law then yields the two ultra-local models:
 \begin{eqnarray}
\label{MFC_LLVC2_11}
\text{longitudinal motion:} & \dot{y}_1 = F_1+ \alpha_1 u_1\\
\label{MFC_LLVC2_22}
\text{lateral motion:} & \ddot{y}_2 = F_2+ \alpha_2 u_2 
\end{eqnarray}
Note the following properties:
\begin{itemize}
\item Equations \eqref{MFC_LLVC2_11}-\eqref{MFC_LLVC2_22} are ``decoupled,''
\item Equation \eqref{MFC_LLVC2_22} is of order $2$ with respect to the derivative of $y_2$. This is the second example of such a property (see Section \ref{A}).
\end{itemize}
Considering the motion in Equation \eqref{MFC_LLVC2_11} (resp. \eqref{MFC_LLVC2_22}), the loop is closed by an iP \eqref{ip} (resp. iPD \eqref{iPD_c}).
\section{Simulation results}
\label{Section_4}
\begin{figure}[!ht]
\centering
\includegraphics[scale=0.5]{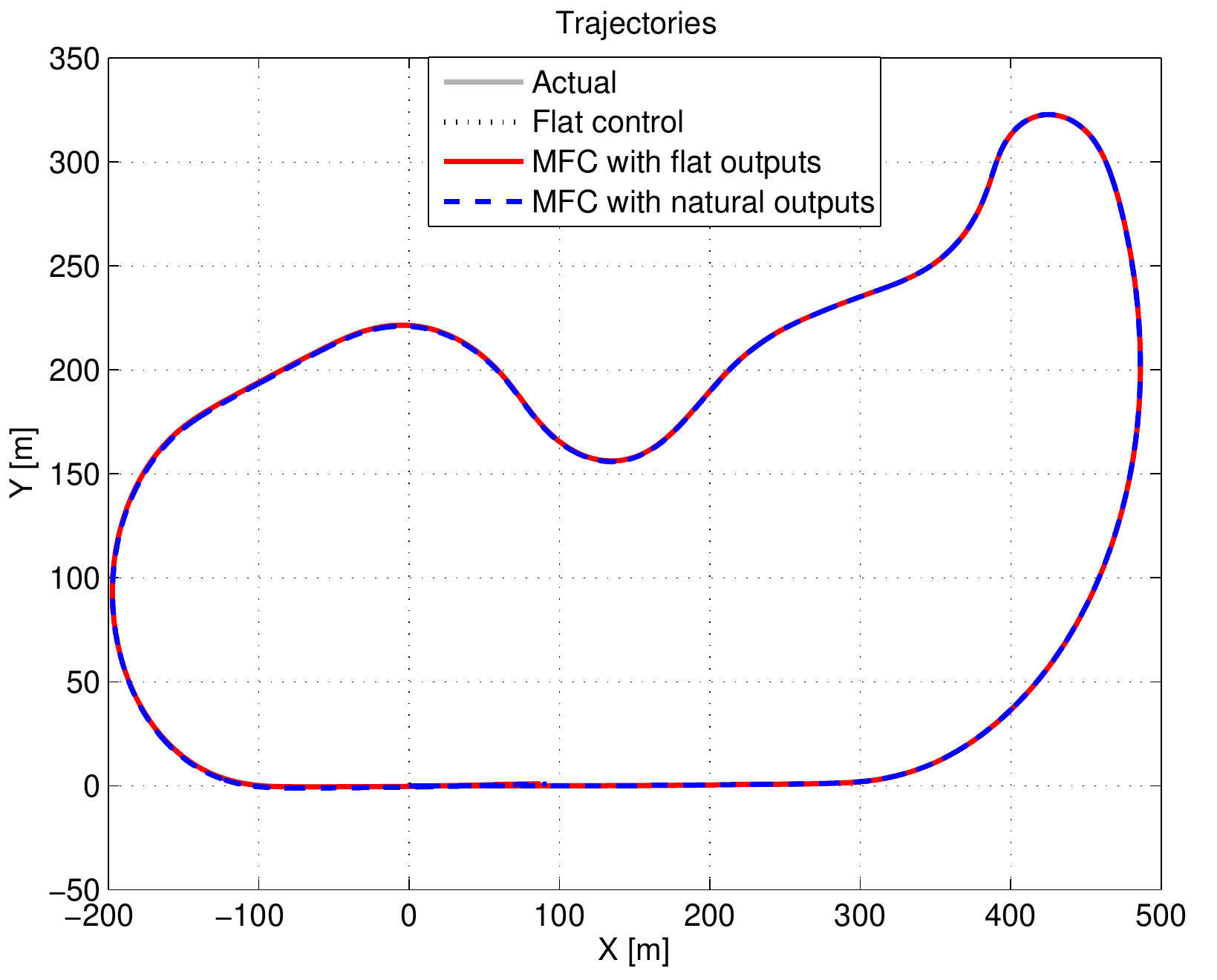}
\vspace{-0.5cm}
\caption{Desired trajectory versus the trajectories of controlled models}
\label{X_Y}
\end{figure}

\begin{figure}[!ht]
\centering
\includegraphics[scale=0.5]{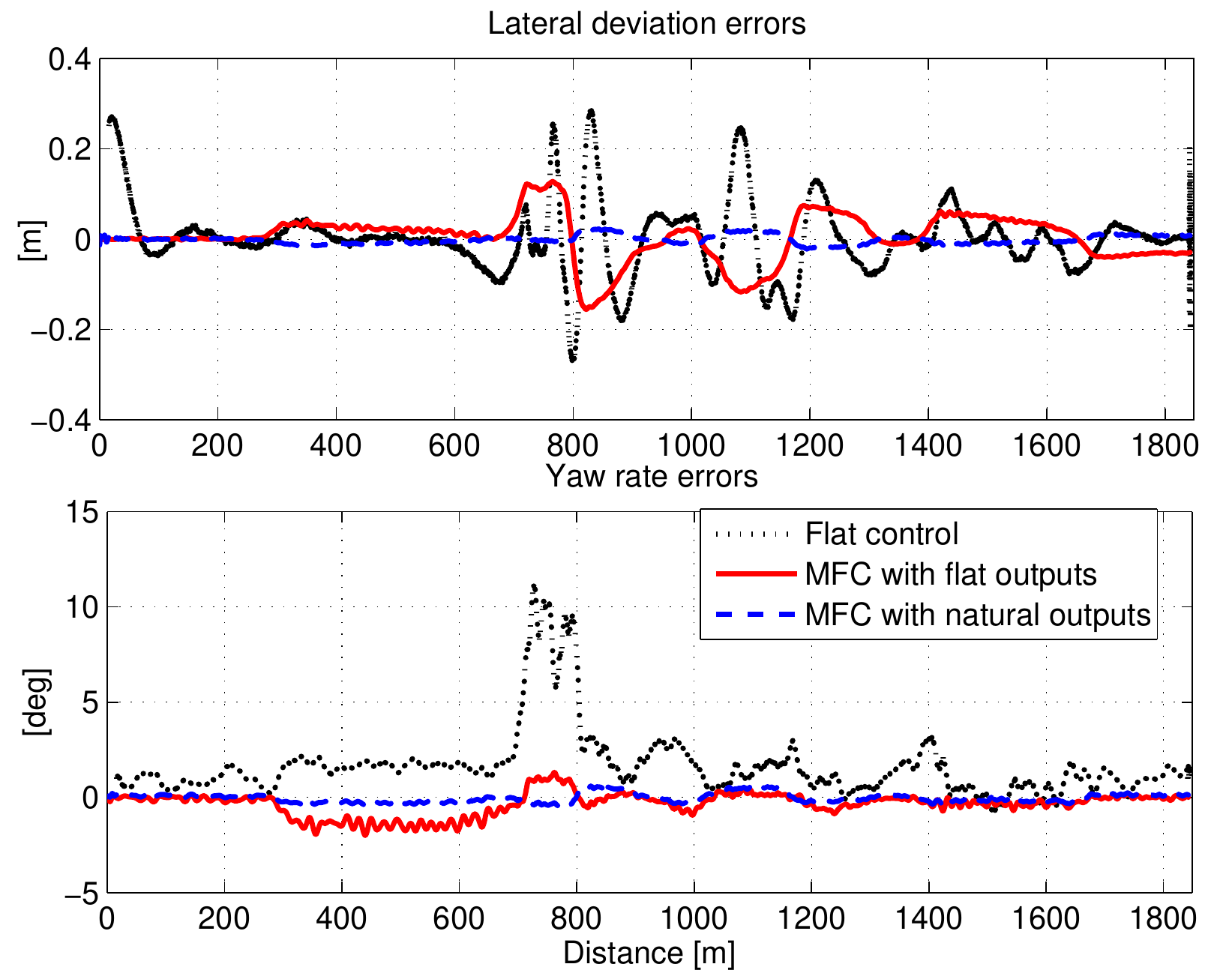}
\vspace{-0.5cm}
\caption{Tacking trajectory errors on lateral deviation and yaw angle}
\label{Ey_Epsi}
\end{figure}

\begin{figure}[!ht]
\centering
\includegraphics[scale=0.5]{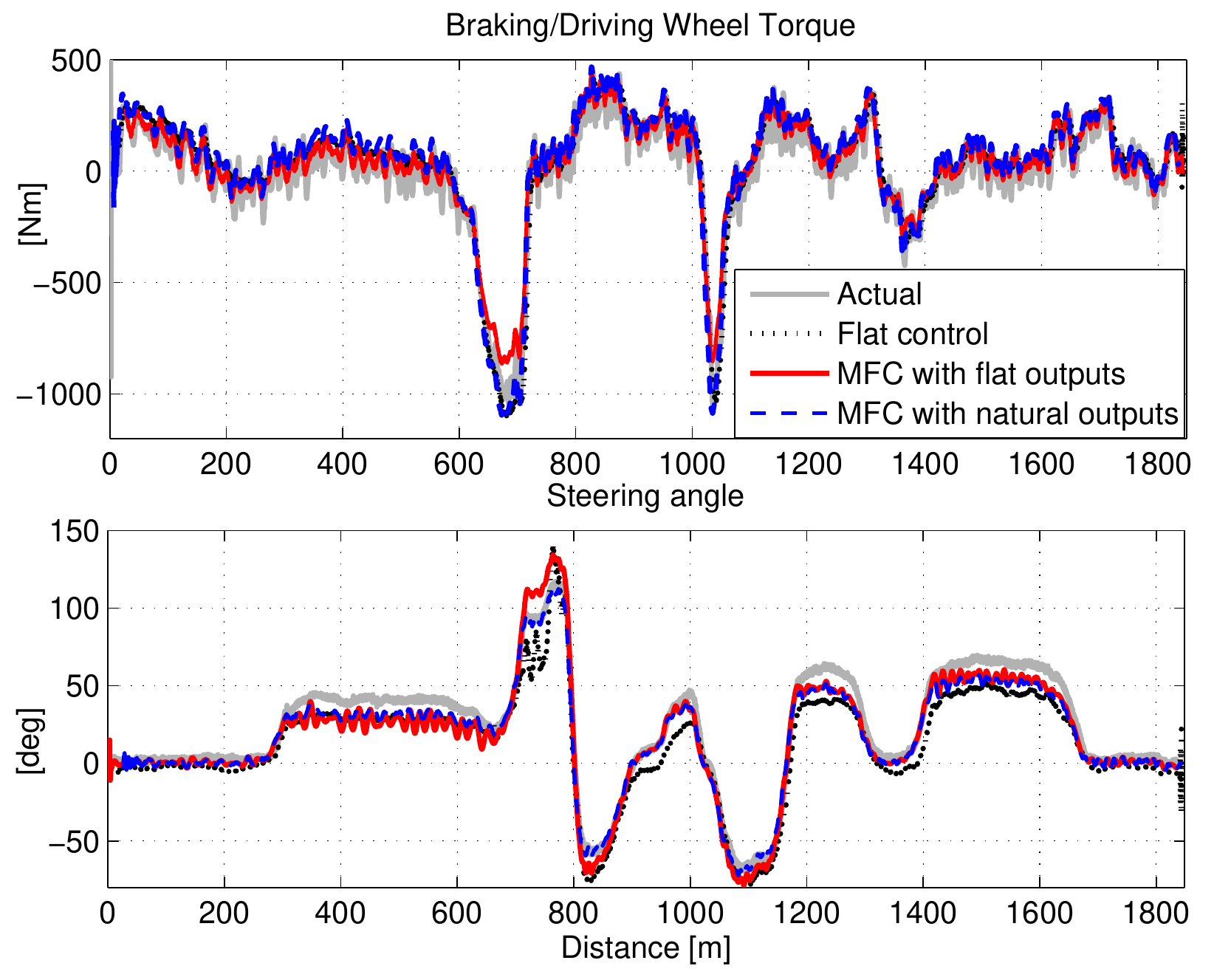}
\vspace{-0.5cm}
\caption{Wheel torques and steering angles control signals: actual and those obtained with new model-free control and the one proposed in \cite{Menhour13b}}
\label{Twheel_Delta}
\end{figure}

The simulations are carried out under MATLAB using a nonlinear 10Dof model of an instrumented Peugeot 406 car as in \cite{Menhour13b}. The reference signals shown by grey curves in Figures \ref{X_Y}, \ref{Ey_Epsi} and \ref{Twheel_Delta} are experimental recorded data. 
Figures \ref{X_Y}, \ref{Ey_Epsi} and \ref{Twheel_Delta} demonstrate that our new model-free setting with natural outputs gives quite satisfactory results. They are much better than those resulting from model-free control associated to the flat outputs. It should be pointed out that the test track which has been considered implies strong lateral and longitudinal dynamical loads. This track involves different types of curvatures linked to straight parts, and all these configurations represent a large set of driving situations. Figure \ref{X_Y} shows that model-free control with natural outputs yields accurate enough behavior for autonomous driving applications. According to the results displayed on Figure \ref{Ey_Epsi}, the lateral error is less than $2 \, cm$. Concerning the yaw angle output, the resulting error is limited to $0.51\, deg$. Note that the lateral errors associated to the flatness-based control and the model-free control with the flat outputs are higher than $10 \, cm$. The resulting yaw angle errors are higher than $10\, deg$ for flatness-based control, and $1\, deg$ for model-free control with flat outputs. Figure \ref{Twheel_Delta} shows that the control signals computed from the model-free control with natural outputs are closed to the actual ones provided by the driver along the track.
\section{Conclusion}
\label{Section_6}
The path followed by the authors for the control of autonomous wheeled cars leads to a design which is both easy to implement and to grasp. The viewpoint, which we now advocate, has been successfully tested under ``normal'' driving conditions through an advanced simulation platform environment (see \cite{ECC2015} for more details). However, in order to be fully convincing, other tests have to be handled on a real embedded vehicle under severe driving conditions. This will be reported in forthcoming publications.

\end{document}